# MONOMIAL IDEALS AND THE SCARF COMPLEX FOR COHERENT SYSTEMS IN RELIABILITY THEORY

BY BEATRICE GIGLIO AND HENRY P. WYNN

*University of Warwick*

A certain type of integer grid, called here an *echelon grid*, is an object found both in coherent systems whose components have a finite or countable number of levels and in algebraic geometry. If $\alpha = (\alpha_1, \ldots, \alpha_d)$ is an integer vector representing the state of a system, then the corresponding algebraic object is a monomial $x_1^{\alpha_1} \cdots x_d^{\alpha_d}$ in the indeterminates $x_1, \ldots, x_d$. The idea is to relate a coherent system to monomial ideals, so that the so-called Scarf complex of the monomial ideal yields an inclusion–exclusion identity for the probability of failure, which uses many fewer terms than the classical identity. Moreover in the "general position" case we obtain via the Scarf complex the tube bounds given by Naiman and Wynn [*J. Inequal. Pure Appl. Math.* (2001) **2** 1–16]. Examples are given for the binary case but the full utility is for general multistate coherent systems and a comprehensive example is given.

**1. Introduction.** The study of network reliability has received increasing attention in the recent decades because of its applications to computer networks and communication systems. Shier [21] points out that one of the most interesting aspects of this subject is the "variety of discrete, combinatorial, and algebraic mathematics that can be found lurking just underneath this practical veneer.... The new approaches developed in the context of network reliability can have ramifications beyond that particular venue."

Naiman and Wynn [16] established a connection between reliability structures and some special topological objects, called abstract tubes. This study led to a second paper [10], where the interactions with computational algebra and Gröbner bases were investigated. This was motivated by the work of Diaconis and Sturmfels [4] that initiated the use of algebraic techniques in probability and statistics. A section on reliability was also included in the monograph ([17], Section 4.4) devoted to establishing or deepening the links existing between algebra and statistics. The work of Dohmen (see, among others, [5, 6, 7]) concentrates more on the combinatorial approach but fruitful and interesting synergies with the abstract tube theory and the algebraic approach are foreseen.

In the present paper the authors investigate a strong link between reliability theory and certain algebraic structures which we believe will eventually unify much of the work just mentioned. The two particular objects covered in this paper are the Scarf complex and the finely graded Hilbert series. These can contribute to finding more efficient computations and identities, similar to those in [16] and [10].

The power of the suggested approach lies in the fact it can be generally applied to any kind of coherent multistate system without any further restrictions. Former approaches proposed in the literature are limited to special cases such as binary systems or source-to-terminal reliability for planar networks. Moreover, because of the link with the tube theory the method immediately gives bounds which are only available, outside the standard Bonferroni case, for certain such cases. Section 5.2 contains a brief discussion of the existing literature and provides examples to illustrate the new methodology.

**2. An overview.** The main point of contact between the disciplines of algebra and reliability theory is an echelon grid of points with nonnegative and integer coordinates. Such a grid (or lattice) can identify a coherent reliability system but also plays a specific role as an algebraic object. An example is given in Figure 1.

A reliability system of $d$ components is a system whose failure or nonfailure state is determined by the state of each of its components. We deal here with the case where each component can assume a discrete and possibly countable number of states $\{0, 1, 2, 3, \ldots\}$ corresponding to increasing levels of efficiency. We refer to such systems as multistate systems. Section 5.5 deals with the case of a continuous distribution on the component performance and more formal definitions about reliability systems are given and discussed in Section 5. A state of the given system is described by a nonnegative integer vector of length $d$. The state space $\mathcal{D}$ ($\mathcal{D} = \mathbb{N}^d$ in this general setting) can therefore be decomposed into the set of configurations that corresponds to failure of the system and the set of configurations that corresponds to nonfailure. We indicate these two sets with $\mathcal{F}$ and $\bar{\mathcal{F}}$, respectively. Obviously $\mathcal{D} = \mathcal{F} \cup \bar{\mathcal{F}}$ (i.e., $\mathcal{F}$ is the complement of $\bar{\mathcal{F}}$ in $\mathcal{D}$).

Coherent systems are systems for which improving the state of a single component cannot lead the system from a nonfailure (or operating) state to a failure state. Equivalently, degrading a component state cannot bring the system from failure to nonfailure. Thus, in the integer representation a coherent system is such that if a state point $\alpha$ belongs to the failure set $\mathcal{F}$, then any other point with component state levels worse (less) than $\alpha$ must belong to $\mathcal{F}$. Similarly the nonfailure set, being the complement of $\mathcal{F}$, must have a similar property: if $\alpha$ is a nonfailure point, then any point with coordinates greater than or equal to $\alpha$'s coordinates must still be a nonfailure point. This is equivalent to saying that the integer representation of the failure set must have an echelon structure; that is, the grid representing $\mathcal{F}$ (or equivalently $\bar{\mathcal{F}}$) must have no "holes."

For example Figure 1 represents a two-component coherent system where the hollow dots correspond to failure states and the filled dots correspond to operating (nonfailure) configurations.

The link with algebra is the construction of the monomial $x^\alpha = x_1^{\alpha_1} x_2^{\alpha_2} \cdots x_d^{\alpha_d}$ with the point $\alpha = (\alpha_1, \alpha_2, \ldots, \alpha_d)$ given by its exponent vector. For example, the filled dots in Figure 1 represent the monomials belonging to a certain monomial ideal $M$. For the introductory material in algebra we have mainly used [3, 12, 14] as references. The echelon property, equivalent to coherence, translates into the algebraic language as the "order ideal" property: if the monomial $x^\alpha$ belongs to a monomial ideal, then so do all the monomials $x^\beta$ such that $x^\alpha$ divides $x^\beta$, that is, such that $\alpha \preceq \beta : \alpha_i \leq \beta_i, i = 1, \ldots, d$. The corner points of the set $\bar{\mathcal{F}}$ (those labeled with their coordinates in Figure 1) are the minimal generators of this monomial ideal but it is interesting to point out that these points play also a specific role as reliability theory objects. They are called *minimal nonfailure points* in Section 5. In this paper we shall only study coherent systems.

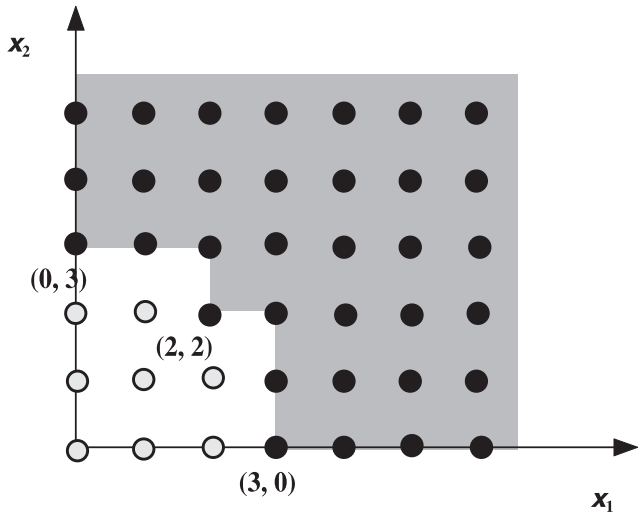

FIG. 1. *Monomials in $M = \langle x_1^3, x_1^2 x_2^2, x_2^3 \rangle$ (black dots) and monomials not in $M$ (hollow dots).*

The main observation, which we will develop in the following sections, is that there are some points on the boundary between the failure set and the nonfailure set that are fundamental to the identification of the points in $\mathcal{F}$ and in $\bar{\mathcal{F}}$. We will see that these special points form a simplicial complex, called the *Scarf complex*, which then plays a key role in both the algebraic geometry and reliability theory.

Sections 3 and 4 give an introduction to the algebraic setting and Section 5 establishes and formalizes the connection with reliability.

### 3. Monomial ideals.

3.1. *Some basics.* Let $\mathbb{K}$ be a field and let $S = \mathbb{K}[x] = \mathbb{K}[x_1, \ldots, x_d]$ be the polynomial ring in $d$ indeterminates. A monomial in $\mathbb{K}[x]$ is a product $x^\alpha = x_1^{\alpha_1} x_2^{\alpha_2} \cdots x_d^{\alpha_d}$ for a vector $\alpha = (\alpha_1, \ldots, \alpha_d) \in \mathbb{N}^d$ of nonnegative integers.

A monomial ideal $M$ is an ideal of $S$ generated by some monomials $x^\alpha$, for $\alpha$ belonging to a subset $A \subseteq \mathbb{N}^d$ (possibly infinite). The ideal $M$ consists therefore of all the polynomials of the form $\sum_{\alpha \in A} h_\alpha x^\alpha$, with $h_\alpha \in S$. In this case we write $M = \langle x^\alpha | \alpha \in A \rangle$.

The following lemma allows us to characterize all the monomials that lie in a given monomial ideal.

LEMMA 1. *Let $M = \langle x^\alpha | \alpha \in A \rangle$ be a monomial ideal. Then a monomial $x^\beta$ lies in $M$ if and only if $x^\beta$ is divisible by $x^\alpha$ for some $\alpha \in A$.*

Observe that the set
$$\alpha + \mathbb{N}^d = \{\alpha + \gamma \mid \gamma \in \mathbb{N}^d\}$$
consists of the exponents of all monomials divisible by $x^\alpha$. This observation allows us to draw pictures of the monomials in a given monomial ideal and visualize this set as a union of positive and integer coordinate points. For example, the filled dots in Figure 1 represent all the monomials in the ideal $M = \langle x_1^3, x_1^2 x_2^2, x_2^3 \rangle$.

A polynomial $f$ is in a monomial ideal $M = \langle x^\alpha \mid \alpha \in A \rangle$ if and only if each term of $f$ is divisible by one of the given generators $x^\alpha$. From this it follows that a monomial ideal is uniquely determined by its monomials; that is, two monomial ideals are the same if and only if they contain the same monomials. For details and proofs see [3, 14].

The main result in this section is given by the following theorem (Dickson's lemma), which states that each monomial ideal is uniquely and finitely generated.

THEOREM 1. *A monomial ideal $M = \langle x^\alpha | \alpha \in A \rangle \subset \mathbb{K}[x_1, \ldots, x_d]$ can be written in the form $M = \langle x^{\alpha_1}, \ldots, x^{\alpha_s} \rangle$, where $\alpha_i \in A$. In particular, $M$ has a finite basis.*

The polynomial ring $S = \mathbb{K}[x]$ can be seen as a $\mathbb{K}$-vector space and therefore it can be decomposed into the direct sum $S = \bigoplus_{\alpha \in \mathbb{N}^d} S_\alpha$, where $S_\alpha$ is the $\mathbb{K}$-span of the monomial $x^\alpha$. Since $S_\alpha \cdot S_\beta \subseteq S_{\alpha+\beta}$, we say that $S$ is an $\mathbb{N}^d$-graded $\mathbb{K}$-algebra. More generally, an $S$-module $M$ is said to be $\mathbb{Z}^d$-graded if $M = \bigoplus_{\beta \in \mathbb{Z}^d} M_\beta$ is a direct sum of $\mathbb{K}$-vector spaces with $S_\alpha \cdot M_\beta \subseteq M_{\alpha+\beta}$. For definitions and properties of modules of polynomial rings see [8].

A monomial ideal $M$ and the corresponding factor ring $S/M$ are both $\mathbb{Z}^d$-graded $S$-modules: $M = \bigoplus_{x^\alpha \in M} S_\alpha$ and $S/M = \bigoplus_{x^\alpha \notin M} S_\alpha$.

The hollow dots in Figure 1 within the nonshaded area form a $\mathbb{K}$-basis for $S/M$. We show in the next section that the Hilbert series $H(S/M; x)$ is the sum of all the monomials that are not contained in $M$.

Given a $\mathbb{Z}^d$-graded module $M$, the $\mathbb{Z}^d$-graded shift $M[\alpha]$ for $\alpha \in \mathbb{Z}^d$ is the $\mathbb{Z}^d$-graded module defined by $M[\alpha]_\beta = M_{\alpha+\beta}$. In particular, the free $S$-module of rank 1 generated in degree $\alpha$ is $S[-\alpha]$. There is an isomorphism between $S[-\alpha]$ and the principal ideal $\langle x^\alpha \rangle \subset S$.

3.2. *Hilbert series and free resolutions.* Given a $\mathbb{Z}^d$-graded module $M_\alpha$ we consider its dimension as a vector space on $\mathbb{K}$ and we call it $\dim_\mathbb{K}(M_\alpha)$. If $\dim_\mathbb{K}(M_\alpha)$ is finite for all $\alpha \in \mathbb{Z}^d$ the (finely graded) *Hilbert series* is defined as the formal power series

$$H(M; x) := H(M; x_1, \ldots, x_d) = \sum_{\alpha \in \mathbb{Z}^d} \dim_\mathbb{K}(M_\alpha) x^\alpha.$$

For example, we have that

$$H(S; x) = \sum_{\alpha \in \mathbb{Z}^d} x^\alpha = \prod_{i=1}^d \frac{1}{1 - x_i}$$

(namely the sum of all monomials in $S$), and for $\alpha \in \mathbb{Z}^d$,

$$\begin{aligned}
H(S[-\alpha]; x) &= \sum_{x^\beta \in \langle x^\alpha \rangle} x^\beta \\
&= x^\alpha (1 + x_1 + x_1^2 + \cdots) \cdots (1 + x_d + x_d^2 + \cdots) \\
&= x^\alpha \frac{1}{1 - x_1} \cdots \frac{1}{1 - x_d} \\
&= \frac{x^\alpha}{\prod_{i=1}^d (1 - x_i)}.
\end{aligned}$$

(1)

For a monomial ideal $M$ we have that

$$H(S/M; x) = \sum_{\alpha \in \mathbb{N}^d \setminus M} x^\alpha,$$

namely the sum of all monomials *not* in $M$.

A *homological complex* of $S$-modules is a sequence

$$\cdots \stackrel{\phi_{i-1}}{\leftarrow} F_{i-1} \stackrel{\phi_i}{\leftarrow} F_i \cdots$$

of $S$-module homomorphisms such that $\phi_{i-1} \circ \phi_i = 0$. A complex is *exact* at the $i$th step if it has no homology there, that is, if $\mathrm{Kernel}(\phi_{i-1}) = \mathrm{Image}(\phi_i)$. The complex is exact if it is exact at the $i$th step for all $i \in \mathbb{Z}$.

A *free resolution* of an $S$-module $M$ is a complex

$$0 \leftarrow F_0 \stackrel{\phi_1}{\leftarrow} F_1 \leftarrow \cdots \leftarrow F_{t-1} \stackrel{\phi_t}{\leftarrow} F_t \leftarrow 0$$

of free modules which is exact everywhere except the 0th step, and such that $M = \mathrm{Coker}(\phi_1) = F_0/\mathrm{Image}(\phi_1)$.

Every $S$-module has a free resolution, with length less than or equal to $d$. If $M$ is $\mathbb{Z}^d$-graded, then it has a $\mathbb{Z}^d$-graded free resolution. If, in addition, $M$ is finitely generated, there is a $\mathbb{Z}^d$-graded resolution $M \leftarrow \mathbb{F}$ in which all the ranks of the $F_i$ are finite and simultaneously minimized. Such an $\mathbb{F}$ is called a *minimal free resolution* of $M$ and is unique up to noncanonical homomorphism.

Given a short exact sequence $0 \leftarrow M'' \leftarrow M \leftarrow M' \leftarrow 0$, the rank nullity theorem of linear algebra implies that $\dim_{\mathbb{K}}(M_\alpha) = \dim_{\mathbb{K}}(M''_\alpha) + \dim_{\mathbb{K}}(M'_\alpha)$ for all $\alpha$, and hence $H(M;x) = H(M'';x) + H(M';x)$. More generally, if

$$0 \leftarrow M \leftarrow F_0 \leftarrow F_1 \leftarrow \cdots$$

is a finite sequence such as a free resolution, then

(2) $$H(M;x) = \sum_i (-1)^i H(F_i;x).$$

In particular, if $M$ is finitely generated, the existence of a finite-rank free resolution for $M$ implies that the Hilbert series of $M$ is a rational function of $x$, because it is an alternating sum of Hilbert series of $S[-\alpha]$ for various $\alpha$. Moreover the denominator can always be taken to be $\prod_i(1-x_i)$.

We show, with an example, the connection between minimal resolution and Hilbert series and why we are particularly interested in it.

EXAMPLE 1. Consider a monomial ideal in two variables: $M = \langle x^{a_1}y^{b_1}, x^{a_2} \times y^{b_2}, \ldots, x^{a_r}y^{b_r}\rangle$ in $S = \mathbb{K}[x,y]$. As mentioned above the finely graded Hilbert series of the factor ring $S/M$ gives the sum of all the monomials not in $M$ (and for the monomial ideal $M$ it gives the sum of all the monomials in $M$). A way to obtain the Hilbert series for $S/M$ is to proceed by inclusion–exclusion, subtracting from the sum of all the monomials in $S$ the list of all the monomials in the quadrants $\langle x^{a_i}y^{b_i}\rangle$ for all $i \in \{1,\ldots,r\}$, then adding back the monomials in the two-way "intersections" $\langle \mathrm{lcm}(x^{a_i}y^{b_i}, x^{a_j}y^{b_j})\rangle = \langle x^{a_i}y^{b_i}\rangle \cap \langle x^{a_j}y^{b_j}\rangle$, then removing the three-way intersections $\langle x^{a_i}y^{b_i}\rangle \cap \langle x^{a_j}y^{b_j}\rangle \cap$ $\langle x^{a_k} y^{b_k} \rangle$, and so on. In this way, after $r$ steps we have counted all the monomials the right number of times. However, for a large number of generators this procedure is far from efficient since many terms cancel out. It is easy to see that in this two-dimensional case it is enough to add back only the principal ideals obtained from couples of adjacent generators.

In the literature the highly nonminimal free resolution of $S/M$ given by the inclusion–exclusion process is called the Taylor resolution, which we describe in more general terms in Section 4.1. The most compact representation yields instead the minimal resolution. Unfortunately a construction of the minimal free resolution has not been found for all arbitrary monomial ideals. In the literature there are only two general constructions for resolving arbitrary monomial ideals: the Taylor resolution and Lyubeznik's subcomplex. We show in Section 4, based on [2, 13, 14], a construction of a minimal free resolution for a specific class of monomial ideals, called *generic monomial ideals*. Also a nonminimal free resolution of $S/M$ based on deformation of the exponents is obtained for nongeneric monomial ideals. Even though this resolution is not minimal it is generally much smaller than Taylor's. We show in Section 5 that in reliability theory the generic monomial ideal has a very natural interpretation but even for binary systems (which typically lead to the nongeneric case) the deformation procedure still gives excellent results.

Notice that the deformation procedure is gratifyingly similar to the perturbation presented in [16] and [10], and the nonminimal resolution obtained via the Scarf complex after deformation leads to the abstract tube formula described in [16] and [10].

**4. Monomial resolutions.** In this section we first describe the construction of a (highly nonminimal) resolution for an arbitrary monomial ideal, called Taylor resolution. Then a method is described to obtain the minimal resolution for generic monomial ideals, and finally the deformation procedure to deal with arbitrary monomial ideals. For more details, see [2, 13, 14].

4.1. *Taylor resolution.* For a monomial ideal $M = \langle m_1, \ldots, m_r \rangle$ and a subset $I \subseteq \{1, \ldots, r\}$ we set $m_I = \text{lcm}(m_i | i \in I)$. Let $\alpha_I \in \mathbb{N}^d$ be the exponent vector of $m_I$ and let $S(-\alpha_I)$ be the free $S$-module with one generator in multidegree $\alpha_I$. Let $F_s$ be the $\mathbb{K}$-vector space whose basis elements $e_I$ correspond to the index sets $I \subseteq \{1, \ldots, r\}$ of length $s$. Define the *differential* $\partial_s : F_s \to F_{s-1}$ by

$$\partial_s(e_I) = \sum_{i \in I} \text{sign}(i, I) \frac{m_I}{m_{I \setminus i}} e_{I \setminus i},$$

where $\text{sign}(i, I)$ is $(-1)^{j+1}$ if $i$ is the $j$th element in the ordering of $I$. It is possible to show (see [8]) that

$$0 \leftarrow F_0 \xleftarrow{\partial_1} F_1 \leftarrow \cdots \leftarrow F_{r-1} \xleftarrow{\partial_r} F_r \leftarrow 0$$

is a free resolution of $S/M$ of length $r$, and it is called *Taylor resolution* of $S/M$. The resolution can be also seen as the $\mathbb{Z}^d$-graded module $\mathbb{F} = \bigoplus_{I \subseteq \{1,\ldots,r\}} S(-\alpha_I)$ that has $2^r$ terms and it is therefore very far from minimal if $r \gg n$.

It is natural to write the Hilbert series for a Taylor resolution of the factor ring $S/M$ as follows:

$$
\begin{aligned}
&H(S/M; x, y) \\
&= H(S; x, y) - H(M; x, y) \\
&= \sum_{\alpha \in \mathbb{N}^d} x^\alpha - \sum_{\substack{x^\alpha \in \langle m_I \rangle \\ I \subseteq \{1,\ldots,r\} \\ |I|=1}} x^\alpha + \sum_{\substack{x^\alpha \in \langle m_I \rangle \\ I \subseteq \{1,\ldots,r\} \\ |I|=2}} x^\alpha + \cdots + (-1)^r \sum_{\substack{x^\alpha \in \langle m_I \rangle \\ I = \{1,\ldots,r\}}} x^\alpha,
\end{aligned}
$$

(3)

where $\langle m_I \rangle$ is the principal ideal generated by the monomial $m_I$. Similarly for the monomial ideal $M$,

(4)
$$
\begin{aligned}
H(M; x, y) \\
= \sum_{\substack{x^\alpha \in \langle m_I \rangle \\ I \subseteq \{1,\ldots,r\} \\ |I|=1}} x^\alpha - \sum_{\substack{x^\alpha \in \langle m_I \rangle \\ I \subseteq \{1,\ldots,r\} \\ |I|=2}} x^\alpha + \cdots + (-1)^{r+1} \sum_{\substack{x^\alpha \in \langle m_I \rangle \\ I = \{1,\ldots,r\}}} x^\alpha.
\end{aligned}
$$

Notice the analogy between the structure of the above formula and the classical inclusion–exclusion formula used in reliability. Simply replace any term $x^\alpha$ by the indicator function of the *orthants* $Q_{\alpha_I} = \{\beta \mid \beta \geq \alpha_I\}$.

We show in the next sections that the Scarf complex allows a more efficient formulation of the above expression and yields exactly the tube formula proposed in [16] and in [10]. For clarity we return to the example of Figure 1.

EXAMPLE 2. We construct now the Taylor resolution and the Hilbert series of $S/M$ where $M = \langle x^3, x^2 y^2, y^3 \rangle$ is the monomial ideal given in Section 3 and represented in Figure 1.

The required resolution is given by

$$0 \leftarrow F_0 \xleftarrow{\partial_1} F_1 \xleftarrow{\partial_2} F_2 \xleftarrow{\partial_3} F_3 \leftarrow 0.$$

The dimension of each vector space $F_s$ ($s = 0, 1, 2, 3$) in the sequence is given by $\binom{d}{s}$ (i.e., the number of subsets $I \subseteq \{1, \ldots, d\}$ of length $s$).

The differentials $\partial_s$ can be expressed via monomial matrices where each column corresponds to an index set $I$ and contains the vector $\partial_s(e_I)$. Therefore we obtain the following:

$$\partial_1: \begin{bmatrix} x^3 & x^2y^2 & y^3 \end{bmatrix};$$
$$\phantom{\partial_1:}\ \ I{=}\{1\}\ \ I{=}\{2\}\ \ I{=}\{3\}$$

$$\partial_2: \begin{matrix} x^3 \\ x^2y^2 \\ y^3 \end{matrix} \left[ \begin{matrix} y^2 & y^3 & 0 \\ -x & 0 & y \\ 0 & -x^3 & -x^2 \end{matrix} \right] ;$$
$$\begin{matrix} \scriptstyle I=\{1,2\} & \scriptstyle I=\{1,3\} & \scriptstyle I=\{2,3\} \\ \scriptstyle m_I=x^3y^2 & \scriptstyle m_I=x^3y^3 & \scriptstyle m_I=x^2y^3 \end{matrix}$$

$$\partial_3: \begin{matrix} x^3y^2 \\ x^3y^3 \\ x^2y^3 \end{matrix} \left[ \begin{matrix} y \\ -1 \\ x \end{matrix} \right] .$$
$$\begin{matrix} \scriptstyle I=\{1,2,3\} \\ \scriptstyle m_I=x^3y^3 \end{matrix}$$

For example, for $I = \{1, 2\}$ it is

$$\partial_2(e_{12}) = \frac{m_{12}}{m_1}e_1 - \frac{m_{12}}{m_2}e_2 = \begin{bmatrix} y^2 \\ -x \\ 0 \end{bmatrix},$$

which corresponds to the first column of the matrix representing $\partial_2$.

Now from (1) and (2) it follows that

$$H(S/M; x, y) = \frac{1 - x^3 - x^2y^2 - y^3 + x^3y^2 + x^3y^3 + x^2y^3 - x^3y^3}{(1-x)(1-y)}$$
$$= 1 + x + y + x^2 + xy + y^2 + x^2y + xy^2,$$

which is exactly the sum of monomials not in $M$, as expected.

Notice that before cancellation of the term $x^3y^3$ the numerator of the Hilbert series for the Taylor resolution contains $2^3 = 8$ terms.

4.2. *Monomial resolution over labeled simplicial complexes and the Scarf complex.* It is possible to define a monomial resolution associated with a simplicial complex labeled by the generators of a given monomial ideal. In this section we show that the resolution of a generic monomial ideal associated with a specific simplicial complex, called the *Scarf complex*, is minimally free. The Scarf complex owes its name to Herbert Scarf, who introduced a similar complex in mathematical economics (see [20]). The importance of this structure for the resolution of a monomial ideal was only recently understood and explained in [2].

DEFINITION 1. Let $V = \{v_1, \ldots, v_r\}$ be a finite set. A (*finite*) *simplicial complex* $\Delta$ on $V$ is a collection of subsets of $V$ such that $I_1 \in \Delta$ whenever $I_1 \subset I_2$ for some $I_2 \in \Delta$, and such that $\{v_i\} \in \Delta$ for $i = 1, \ldots, r$. The elements of $\Delta$ are called *faces*, and the *dimension* $\dim(I)$ of a face $I$ is the number $|I| - 1$. The *dimension* of the simplicial complex $\Delta$ is $\dim(\Delta) = \max\{\dim(I) : I \in \Delta\}$.

Note that the empty set $\varnothing$ is a face of dimension $-1$ of any nonempty simplicial complex. Faces of dimension 0 and 1 are called *vertices* and *edges*, respectively. The maximal faces under inclusion are called *facets* of the simplicial complex.

For a monomial ideal $M$ in $S = \mathbb{K}[x_1, \ldots, x_d]$ we can construct a special simplicial complex whose vertices are labeled with the generators of $M$. Additionally we label each face $I$ of the simplicial complex by the least common multiple $m_I = \mathrm{lcm}(m_i | i \in I)$ of its vertices and we restrict to the faces with unique labels. This leads to the following definition.

DEFINITION 2. For a monomial ideal $M = \langle m_1, \ldots, m_r \rangle$ we call the *Scarf complex* the simplicial complex $\Delta_M$ consisting of sets of minimal generators with unique labels

$$\Delta_M = \{I \subseteq \{1, \ldots, r\} \mid m_I \neq m_J \text{ for all } J \subseteq \{1, \ldots, r\} \text{ other than } I\}.$$

It is not difficult to see that a Scarf complex in $S = \mathbb{K}[x_1, \ldots, x_d]$ is still a simplicial complex and has dimension at most $d - 1$.

For $i \in \mathbb{Z}$ let $F_i(\Delta_M)$ be the set of $i$-dimensional faces of $\Delta_M$ and let $\mathbb{K}^{F_i(\Delta_M)}$ be a $\mathbb{K}$-vector space whose basis elements $e_I$ correspond to the $i$-faces $I \in F_i(\Delta_M)$. Let $\mathbb{F}_{\Delta_M}$ be the so-called $\mathbb{N}^d$-graded chain complex of $\Delta_M$ over $S$ obtained as

$$0 \leftarrow \mathbb{K}^{F_{-1}(\Delta_M)} \stackrel{\partial_0}{\leftarrow} \cdots \leftarrow \mathbb{K}^{F_{i-1}(\Delta_M)} \stackrel{\partial_i}{\leftarrow} \mathbb{K}^{F_i(\Delta_M)} \leftarrow \cdots \stackrel{\partial_{n-1}}{\leftarrow} \mathbb{K}^{F_{n-1}(\Delta_M)} \leftarrow 0$$

with differentials as in Section 4.1,

$$\partial_s(e_I) = \sum_{i \in I} \mathrm{sign}(i, I) \frac{m_I}{m_{I \setminus i}} e_{I \setminus i}.$$

For a special class of monomial ideals called *generic*, the complex defined by the simplicial complex $\Delta_M$ is minimally free as Theorem 2 states (for the proof see [2]).

DEFINITION 3. A monomial ideal $M$ is called *generic* if no variable $x_i$ appears with the same nonzero exponents in two distinct minimal generators of $M$.

Note that in [13, 14] the authors prove that it is possible to define genericity in a weaker way and still achieve the same results that we are going to present. They therefore refer to the definition we give above as "strong genericity." For the purpose of this article the stronger version is not only sufficient but also easier and more appropriate to use for application to reliability. The requirement that the monomial ideal is generic might seem quite strong. In practice though almost all monomial ideals are generic, in the sense that those which fail to be generic lie on finitely many hyperplanes in the space of exponents. From the reliability point of view the requirement is not too strong in the continuous distribution case (see Section 5.5), while in the binary case the deformation procedure can be used to obtain a resolution (typically nonminimal) as we shall see from examples.

THEOREM 2. *For a generic monomial ideal $M$ the complex $\mathbb{F}_{\Delta_M}$ defined by the Scarf complex $\Delta_M$ is the minimal free resolution of $S/M$ over $S$. The $\mathbb{N}^d$-graded Hilbert series of $S/M$ (i.e., the sum of all monomials not in $M$) is*

$$
(5) \qquad \frac{\sum_{I \in \Delta_M} (-1)^{|I|} \, m_I}{(1-x_1)\cdots(1-x_d)}
$$

*and there are no cancellations in the alternating sum in the numerator.*

Theorem 2 allows us to write the Hilbert series for the factor ring $S/M$ (and for the monomial ideal $M$) in a much more parsimonious way than in the Taylor resolution. The best way to see the difference is to compare (3) and (4) with the following expressions for $H(M;x,y)$ and $H(S/M;x,y)$:

$$
(6) \quad \begin{aligned} H(S/M;x,y) &= \sum_{\alpha \in \mathbb{N}^d} x^\alpha - \sum_{\substack{x^\alpha \in \langle m_I \rangle \\ I \in \Delta_M \\ |I|=1}} x^\alpha \\ &\quad + \sum_{\substack{x^\alpha \in \langle m_I \rangle \\ I \in \Delta_M \\ |I|=2}} x^\alpha + \cdots + (-1)^r \sum_{\substack{x^\alpha \in \langle m_I \rangle \\ I = \{1,\ldots,r\}}} x^\alpha; \end{aligned}
$$

$$
(7) \quad H(M;x,y) = \sum_{\substack{x^\alpha \in \langle m_I \rangle \\ I \in \Delta_M \\ |I|=1}} x^\alpha - \sum_{\substack{x^\alpha \in \langle m_I \rangle \\ I \in \Delta_M \\ |I|=2}} x^\alpha + \cdots + (-1)^{r+1} \sum_{\substack{x^\alpha \in \langle m_I \rangle \\ I = \{1,\ldots,r\}}} x^\alpha.
$$

The Taylor resolution can be obtained from the full simplicial complex (with cardinality $2^r$) given by the set of *all* subsets of $\{1,\ldots,r\}$, while in the Scarf resolution the summation is done only on the subsets belonging to $\Delta_M$. Let us see now with an example how we can obtain the minimal resolution via the Scarf complex for the monomial ideal considered in Example 2 (Figure 1). The example shows a "by hand" construction of the Scarf complex while for most of the other examples the authors made continual use of available functions in CoCoA. The algorithms behind these functions are related to the construction of the Hilbert series (see [14]) and seem reasonably efficient. Further research is needed to establish theoretical or experimental results on formal computational complexity or actual run time. We do not develop the computational aspects in this paper.

EXAMPLE 3. The Scarf complex for the monomial ideal $M = \langle x^3, x^2y^2, y^3 \rangle$ can be obtained in the following way:

$$I = \{1\} \to m_I = x^3;$$
$$I = \{2\} \to m_I = x^2y^2;$$
$$I = \{3\} \to m_I = y^3;$$
$$I = \{1,2\} \to m_I = x^3y^2;$$
$$I = \{1,3\} \to m_I = x^3y^3;$$
$$I = \{2,3\} \to m_I = x^2y^3;$$
$$I = \{1,2,3\} \to m_I = x^3y^3.$$

Thus the Scarf complex $\Delta_M$ is given by

$$\Delta_M = \{\{1\}, \{2\}, \{3\}, \{1,2\}, \{2,3\}\}.$$

The index subsets that do not appear in $\Delta_M$ are those corresponding to the monomial $x^3y^3$, and notice that this is exactly the term that cancels in Example 2.

The resolution is given by

$$0 \leftarrow S \xleftarrow{[x^3, x^2y^2, y^3]} S^3 \xleftarrow{\begin{vmatrix} y^2 & 0 \\ -x & y \\ 0 & -x^2 \end{vmatrix}} S^2 \leftarrow 0,$$

and the Hilbert series (which can be read off directly from the labels of the Scarf complex) equals

$$H(S/M; x) = \frac{1 - x^3 - x^2y^2 - y^3 + x^3y^2 + x^2y^3}{(1-x)(1-y)}.$$

Thus the number of terms in the numerator of the Hilbert series is six. The gain in comparison with the Taylor resolution, which returned eight terms, is not large for this small example but it becomes huge in high-dimensional problems.

4.3. *Deformation of exponents.* As pointed out in Section 4.2, when the generators of the monomial ideal are in generic position the Scarf complex leads to the minimal resolution. For arbitrary monomial ideals the construction described in this section can be used to produce a (typically nonminimal) resolution via deformation of the exponent vectors of the generators. There is no guarantee that the resolution so obtained is minimal, but the authors' experience is that there is still a very great gain in resolution over the original Taylor complex. The use of the deformation procedure combined with the use of the CoCoA macros mentioned in Section 4.2 seems very powerful and certainly superior to complete enumeration followed by cancellation. The example presented in Section 5.4 shows a typical result where, even after deformation, a considerable reduction in the number of terms in the inclusion–exclusion formula is achieved (from 511 to 31). Another advantage of the approach presented here is that, even in the cases where minimality is not guaranteed, the method, unlike other methods described in the literature, still provides improved reliability bounds, as described in Section 5.2.

For an arbitrary monomial ideal $M = \langle m_1, \ldots, m_r \rangle$, let $\{\alpha_i = (\alpha_{i1}, \ldots, \alpha_{id}) | 1 \leq i \leq r\}$ be the exponent vectors of the minimal generators of $M$. Choose vectors $\varepsilon_i = (\varepsilon_{i1}, \ldots, \varepsilon_{id}) \in \mathbb{R}^d$ for $1 \leq i \leq r$ such that, for all $i$ and $s \neq t$, the numbers $a_{is} + \varepsilon_{is}$ and $a_{it} + \varepsilon_{it}$ are distinct and $a_{is} + \varepsilon_{is} < a_{it} + \varepsilon_{it}$ implies $a_{is} \leq a_{it}$. Each vector $\varepsilon_i$ defines a monomial $x^{\varepsilon_i} = x_1^{\varepsilon_{i1}} \cdots x_d^{\varepsilon_{id}}$ with real exponents. We define the generic monomial ideal $M_\varepsilon$ in a polynomial ring with real exponents as follows:

$$M_\varepsilon = \langle m_1 x^{\varepsilon_1}, m_2 x^{\varepsilon_2}, \ldots, m_r x^{\varepsilon_r} \rangle.$$

We call $M_\varepsilon$ a generic deformation of $M$. Let $\Delta_{M_\varepsilon}$ be the Scarf complex of $M_\varepsilon$. We now label the vertex of $\Delta_{M_\varepsilon}$ corresponding to $m_i x^{\varepsilon_i}$ with the original monomial $m_i$. Let $\mathcal{F}_\varepsilon$ be the complex of $S$-modules defined by this labeling of $\Delta_{M_\varepsilon}$. Then the following result can be obtained (see the proof in [2]).

THEOREM 3. *The complex $\mathbb{F}_\varepsilon$ is a free resolution of $S/M$ over $S$.*

There is a simple way to deform $M$. This is the deformation given in [2] and is equivalent to the perturbation proposed in [10]. The method consists in picking an integer $v > r$ and deforming $M$ using $\varepsilon_{ij} = i/v$.

## 5. Algebra and reliability.

5.1. *Coherent systems.* We consider here systems as defined in [17] and [1] but with some slight differences in notation.

DEFINITION 4. A *system* is a set $S$ of $d$ components. We code their increasing efficacy levels with the integers $\{0, 1, 2, \ldots\}$.

An *outcome* is a nonnegative integer vector of length $d$ describing the state of each component. We call $\mathcal{D}$ the set of all possible outcomes.

A *failure outcome* is an outcome which leads to failure of the system $S$.

The *failure set* $\mathcal{F}$ is the set of all failure outcomes. The *nonfailure set* is the complement of $\mathcal{F}$ in $\mathcal{D}$: $\bar{\mathcal{F}} = \mathcal{D} \setminus \mathcal{F}$.

In most cases it is natural to assume that replacing a component by a component at a higher efficacy level will not lead to a deterioration of the system. Systems for which such an assumption is valid are called *coherent systems*. We consider only this type of system. For example, if the outcome $(1, 2, 1, 3, 0)$ is an operating (nonfailure) outcome then also $(1, 2, 2, 4, 0)$ must be an operating configuration. In terms of the usual partial order relation $\preceq$ in $\mathbb{Z}^d$, for a coherent system we write that if $\alpha \in \bar{\mathcal{F}}$, then $\beta \in \bar{\mathcal{F}}$ for all $\beta$ such that $\alpha \preceq \beta$ (i.e., $\alpha_i \leq \beta_i$, $i = 1, \ldots, d$). In the same way, in a failure configuration, replacing an operating component with a failed one cannot improve the system. Therefore, if $(0, 1, 0, 0, 0)$ is a failure outcome then $(0, 0, 0, 0, 0)$ must also be. This observation leads to the following definition.

DEFINITION 5. The minimal failure points are maximal points in the partial order relationship defined by $\preceq$:

$$\mathcal{F}^* = \{\alpha \in \mathcal{F} | \nexists \beta \in \mathcal{F} : \alpha \preceq \beta\}.$$

The use of the word "minimal" in the definition is justified by the fact that in the binary network literature the minimal failure points (sometimes called *minimal cuts*) are the minimal set of components whose failure ensures the failure of the system. The collection of all the minimal failure points is called a *minimal failure set* and is indicated with the symbol $\mathcal{F}^*$.

Similarly the minimal points (according to $\preceq$) in the nonfailure set are called *minimal nonfailure points* ( *paths* in network theory):

$$\bar{\mathcal{F}}^* = \{\alpha \in \bar{\mathcal{F}} | \nexists \beta \in \bar{\mathcal{F}} : \beta \preceq \alpha\}.$$

As mentioned in the Introduction it is easy to see how the failure set for a coherent system can be represented by an integer grid. Furthermore this grid has an echelon structure since the system is coherent. A similar observation, but with reverse inequalities, can be made for the nonfailure set.

The parallel between a monomial ideal and a nonfailure set, and between minimal nonfailure points and minimal generators, is now fully established.

PROPOSITION 1. *Given a system $S$ of $d$ components with nonfailure set $\bar{\mathcal{F}}$ and minimal nonfailure set $\bar{\mathcal{F}}^*$:*

1. *The minimal nonfailure points in $\bar{\mathcal{F}}^*$, seen as the exponent vectors of monomials in $\mathbb{K}[x_1, \ldots, x_d]$, are the (exponents of the) minimal generators of a monomial ideal $M$.*
2. *The points in $\bar{\mathcal{F}}$ represent the (exponents of the) monomials belonging to the monomial ideal $M$ generated by the minimal failure points.*
3. *The points in $\mathcal{F}$ represent the (exponents of the) monomials belonging to the factor ring $S/M$.*

Thus we allow ourselves to use interchangeably the expressions "nonfailure set $\bar{\mathcal{F}}$" and "monomial ideal generated by $\bar{\mathcal{F}}^*$," or "minimal failure points" and "minimal generators," even though, strictly speaking, we are dealing with monomials in one case and their exponents in the other one.

The main interest in reliability theory is to give a measure of the performance of a system $S$ by evaluating the probability of the failure (or nonfailure) of the system.

DEFINITION 6. For a system $S$ the reliability function $\mathcal{R}(S)$ is the probability of the nonfailure set $\bar{\mathcal{F}}$; the unreliability $\mathcal{U}(S)$ is the probability of the failure set $\mathcal{F}$. Clearly, we only need to concentrate on one of the two performance measures, since $\mathcal{U}(S) = 1 - \mathcal{R}(S)$.

In the literature it is possible to find many techniques to obtain the reliability of a system, but the most common problem encountered is the computational effort required to evaluate this probability for large size problems.

One of the methods proposed in the literature to calculate the reliability function in the case of a finite state space is due to Moore and Shannon [15] and is based on state space enumeration. Thus the reliability function for a network $G$ with nonfailure set $\bar{\mathcal{F}}$ can be written as

$$\mathcal{R}(S) = \sum_{\omega \in \mathcal{D}} I_{\bar{\mathcal{F}}}(\omega) \Pr(\omega),$$

where $\mathcal{D}$ is the full set of network states and $I_{\bar{\mathcal{F}}}(\omega)$ takes the value 1 when $\omega$ belongs to $\bar{\mathcal{F}}$, 0 otherwise. The approach is impractical because the space $\mathcal{D}$ has cardinality equal to $m^d$, if $d$ is the number of components and $m$ the number of levels that each component can assume.

The reliability function can be formulated in terms of the minimal nonfailure points. For a minimal nonfailure point $\alpha \in \bar{\mathcal{F}}^*$ we indicate with $Q_\alpha$ the event "the system is in one of the states $\beta$, with $\beta \succeq \alpha$." This can also be described as the event that the system state is a point in the *orthant* that has the point $\alpha$ as corner point,

$$Q_\alpha = \{\beta | \beta \succeq \alpha\}.$$

Notice that the orthants $Q_\alpha$ correspond exactly to the set of monomials belonging to the principal ideals $\langle x^\alpha \rangle$. Thus the nonfailure set can be written as a union of the orthants $Q_\alpha$ based on the minimal nonfailure points $\alpha \in \bar{\mathcal{F}}^*$,

$$\bar{\mathcal{F}} = \bigcup_{\alpha \in \bar{\mathcal{F}}^*} Q_\alpha.$$

The reliability function can therefore be obtained as the probability of a union of orthants. Since this union is not disjoint the classical approach to this problem is to use the inclusion–exclusion formula,

$$\begin{aligned}
\mathcal{R}(S) &= \mathrm{Prob}\left( \bigcup_{\alpha \in \bar{\mathcal{F}}^*} Q_\alpha \right) \\
&= \sum_{\alpha \in \bar{\mathcal{F}}^*} \mathrm{Prob}(Q_\alpha) - \sum_{\alpha, \alpha' \in \bar{\mathcal{F}}^*} \mathrm{Prob}(Q_\alpha \cap Q_{\alpha'}) + \cdots \\
&\quad + (-1)^{r+1} \mathrm{Prob}(Q_\alpha \cap Q_{\alpha'} \cap \cdots),
\end{aligned}$$

where $r = |\bar{\mathcal{F}}^*|$ is the cardinality of the minimal nonfailure set.

5.2. *The method.* To obtain the improved version of the inclusion–exclusion formula for a system $S$ with minimal nonfailure set $\bar{\mathcal{F}}^* = \{\alpha_1, \ldots, \alpha_r\}$ we only need the Scarf complex $\Delta_M$ associated with the monomial ideal $M$ generated by the monomials $x^{\alpha_1}, \ldots, x^{\alpha_r}$. The minimal resolution and the Hilbert series provide the background for our calculations. The Scarf complex shows directly which least common multiples are needed to be included in the identity. We decompose the complex in terms of the dimensions of its faces and we write

$$\Delta = \{\Delta_0, \Delta_1, \ldots, \Delta_{d-1}\} \qquad \text{with } d \leq n,$$

where $\Delta_s$ is the set of faces of $\Delta$ of dimension $s+1$ (thus $\Delta_0$ is the set of vertices, $\Delta_1$ is the set of edges, etc.).

From (7) we obtain the improved version of the inclusion–exclusion formula,

$$\begin{aligned}
\mathcal{R}(S) &= \mathrm{Prob}\left( \bigcup_{\alpha \in \bar{\mathcal{F}}^*} Q_\alpha \right) \\
&= \sum_{m_I \in \Delta_0} \mathrm{Prob}(Q_{\alpha_I}) - \sum_{m_I \in \Delta_1} \mathrm{Prob}(Q_{\alpha_I}) + \cdots \\
&\quad + (-1)^{d-1} \sum_{m_I \in \Delta_{d-1}} \mathrm{Prob}(Q_{\alpha_I}),
\end{aligned} \tag{8}$$

where $m_I = \mathrm{lcm}(m_i | i \in I)$ is as defined in Section 4.2 and $\alpha_I$ is the exponent vector of the monomial $m_I$.

This formula is equivalent to the improved inclusion–exclusion formula obtained in [16] (tube identity) for orthant arrangements. The authors prove there that truncating the formula at even and odd level leads to upper and lower bounds for the reliability,

$$\sum_{j=0}^{r+1}(-1)^j \sum_{m_I \in \Delta_j} \text{Prob}(Q_{\alpha_I}) \leq \mathcal{R}(S) \leq \sum_{j=0}^{r}(-1)^j \sum_{m_I \in \Delta_j} \text{Prob}(Q_{\alpha_I}),$$
(9)
$$0 \leq r \leq d-1, r \text{ even}.$$

In [16] it is proved that such bounds are always at least as tight as the bounds from truncating the usual inclusion–exclusion lemma. It is of some interest that these inequalities have not, to the authors' knowledge, been established in the algebraic literature.

5.3. *The binary case.* The methods of this paper were introduced to be applied to any general multistate coherent system and in the next section we give an example. We discuss first, however, the standard binary case on which there is much literature.

This literature divides into three broad classes: cases derived from networks, special nonnetwork cases (such as $k$-out-of-$n$) and arbitrary binary examples. For network and some special nonnetwork cases there are identities and bounds competitive with those presented in this paper. The best published results contain identities with the maximum amount of cancellation in the classical inclusion–exclusion identities. An example is the work by Satyanarayana and Prabhakar (see [18, 19]). For bounds the recent results of Dohmen [5, 6, 7] also using tube identities are based on a somewhat different construction and are specialized to certain classes of networks.

As mentioned, the immediate problem with using the Scarf resolution for binary systems is that the genericity condition only holds in trivial cases and the deformation procedure has to be used. The algebraic theory tells us that this does not guarantee a minimal resolution and, as pointed out in [10], different deformations give rise to different Scarf complexes. Thus, some deformations can be better than others, in the sense that they lead to a smaller Scarf complex and therefore to a more concise reliability identity. The authors are currently investigating the possibility of identifying the "optimal" deformation, that is, the deformation that leads to the minimal resolution. As the general theoretical minimal resolution problem appears to be unsolved in algebra an early theoretical solution is not expected, but fast computational methods look promising.

Notice that in the binary case if we assume independence among the component failures and we indicate by $p_i$ the probability of nonfailure of component $i$, the probability of an orthant $Q_{\alpha_I}$ associated with the index set $I$ in the Scarf complex is easily obtained as

$$\text{Prob}(Q_{\alpha_I}) = \prod_{\alpha_{Ii}=1} p_i. \tag{10}$$

To summarize, the strategy for the binary case which derives from this paper is as follows: (i) describe the failure set; (ii) perturb the minimal failure points to general position; (iii) derive the Scarf complex; (iv) consider the labels (monomials) associated with the Scarf complex; (v) obtain the orthant probabilities; (vi) derive the final identities or bounds.

EXAMPLE 4 (Binary network). We illustrate our method by considering an eight-component binary system defined by the following nine minimal nonfailure points:

$$\bar{\mathcal{F}}^* = \{(1,0,0,0,0,1,0,0),(1,0,0,1,0,0,1,0),(0,1,0,1,0,1,0,0),$$
$$(1,0,0,1,1,0,0,1),(0,1,0,0,0,0,1,0),(0,0,1,1,1,1,0,0),$$
$$(0,1,0,0,1,0,0,1),(0,0,1,0,1,0,1,0),(0,0,1,0,0,0,0,1)\}.$$

This example is taken from [6, 7]. The corresponding monomial ideal in $\mathbb{K}[x_1,\ldots,x_8]$ is generated by the monomials corresponding to the minimal nonfailure points:

$$M = \langle x_1x_6, x_1x_4x_7, x_2x_4x_6, x_1x_4x_5x_8, x_2x_7, x_3x_4x_5x_6, x_2x_5x_8, x_3x_5x_7, x_3x_8\rangle.$$

After perturbation (with $\varepsilon_i = i/10$) and ranking, we obtain the following generic monomial ideal:

$$\widetilde{M} = \langle x_1^8 x_2^5 x_3^5 x_4^4 x_5^4 x_6^8 x_7^5 x_8^5,\ x_1^7 x_2^4 x_3^4 x_4^8 x_5^3 x_6^5 x_7^8 x_8^4,\ x_1^5 x_2^8 x_3^3 x_4^7 x_5^2 x_6^7 x_7^4 x_8^3,$$
$$x_1^6 x_2^3 x_3^2 x_4^6 x_5^8 x_6^4 x_7^3 x_8^8,\ x_1^4 x_2^7 x_3 x_4^3 x_5^3 x_6^3 x_7^7 x_8^2,\ x_1^3 x_2^2 x_3^8 x_4^5 x_5^7 x_6^6 x_7^2 x_8,$$
$$x_1^2 x_2^6 x_4^2 x_5^6 x_6^2 x_7 x_8^7,\ x_1 x_2 x_3^7 x_4 x_5^5 x_6 x_7^6,\ x_3^6 x_8^6 \rangle.$$

The corresponding Scarf complex $\Delta_{\widetilde{M}}$ contains 103 elements, which are the index subsets of the following six facets:

$$\{\{12479\},\{13689\},\{12579\},\{12589\},\{13579\},\{13589\}\},$$

where $1,2,\ldots,9$ are the labels of the minimal generators of $M$.

The reliability function is obtained following (8) and (10) and it is based on the probabilities of the orthants identified by the Scarf complex. For example, the index set $I = \{127\}$ corresponds to the least common multiple of the first, second and seventh minimal generators: $m_I = x_1x_2x_4x_5x_6x_7x_8 = \text{lcm}(x_1x_6, x_1x_4x_7, x_2x_5x_8)$. In the independence case, using (10), the probability of the corresponding orthant is expressed in terms of the probability of nonfailure of the single components, $p_i, i = 1,\ldots,8$,

$$\alpha_I = (1,1,0,1,1,1,1,1), \qquad \text{Prob}(Q_{\alpha_I}) = p_1p_2p_4p_5p_6p_7p_8,$$

where $\alpha_I$ is the exponent of the monomial $m_I$.

Notice that, interestingly, the Scarf complex here coincides exactly with the tube simplicial complex obtained by Dohmen [6, 7] in the case of a binary planar network whose reliability function is defined as source-to-terminal reliability. However, the Scarf complex method can be applied to much more generally defined coherent systems, as for example the binary nonnetwork case or the multistate case.

EXAMPLE 5 (Binary nonnetworks). The recent literature contains several efficient algorithms to calculate the reliability of special classes of coherent systems. For example, the method presented in [19] gives an efficient algorithm that can be used to evaluate the source-to-terminal reliability of binary networks. The work by Dohmen, mentioned in the previous example, provides reliability bounds and the exact formula for a wider class of networks. However, there are other binary structures that are not networks. A $k$-out-of-$n$ system is well studied but it is only one of the many possible examples of a nonnetwork system. In our setting a $k$-out-of-$n$ system can be represented, as any other coherent system, by a monomial ideal. For instance, the ideal

$$M = \langle x_1x_2, x_1x_3, x_1x_4, x_2x_3, x_2x_4, x_3x_4 \rangle$$

is the monomial ideal corresponding to the nonfailure set of a 2-out-of-4 system. In general, different situations can arise where the minimal nonfailure points do not correspond to the paths of any binary networks. For example, it can be seen that the minimal nonfailure set $\bar{\mathcal{F}}^* = \{(1,1,0,0,0),(0,1,1,0,0),(0,0,1,1,0),(1,0,0,1,1),(0,1$
cannot derive from any network. As in the previous example, labeling the minimal nonfailure points $(1,\ldots,5)$, after deformation and ranking the resulting Scarf complex is obtained:

$$\Delta = \{\{1235\}, \{123\}, \{135\}, \{125\}, \{235\}, \{345\},$$
$$\{13\}, \{12\}, \{15\}, \{34\}, \{23\}, \{25\}, \{35\}, \{45\}, \{1\}, \{2\}, \{3\}, \{4\}, \{5\}\}.$$

5.4. *The multistate case.* As explained, the full power of the method is exhibited for multistate systems where the literature on identities is sparser than for the binary case and bounds even more so. The same steps (i)–(iv) are used as in the binary case. For multistate systems the degree of nongenericity tends to be relatively lower because the multiple levels tend to separate out the points better. In the next section we show that with some care the generic setting can also be extended to cover the continuous case.

In Section 5 we introduced the notion of a system without imposing any limit on the number of state levels that each component can assume. The orthants $Q_\alpha$ as defined above were infinite portions of the state space. However, in the situation where each component has only a finite number of states the state space $\mathcal{D}$ is finite and the nonfailure set $\bar{\mathcal{F}}$ is given by the intersection of the monomials in the monomial ideal generated by the set $\bar{\mathcal{F}}^*$ with the state space $\mathcal{D}$. The formula in (8) is still valid if we take care to intersect each orthant with the set $\mathcal{D}$. Full details are given in [11], where the finite case is described in a more formal setting.

We consider here a multistate system where the nonfailure set is defined by setting a cutoff point $c$ for an increasing function $\Psi$ defined on the state space $\mathcal{D}$,

$$\bar{\mathcal{F}}_c = \{\alpha \in \mathcal{D} | \Psi(\alpha) \geq c\}.$$

An example of such a setting is given by a situation where the profit $\Psi$ deriving from a certain system is a function of the performance of each single component. In addition, we consider the case where the performance (profit) of the system is mostly affected by the performance of certain components or combination of components. For example, consider the following profit $\Psi$ for a given four-component system:

(11) $$\Psi = \alpha_1 + \alpha_2 + 4\alpha_3 + 5\alpha_4 + 2\alpha_3\alpha_4,$$

where the performance $\alpha_i$ of the $i$th component can take values in $\{0, 1, 2, 3\}$. The coefficient 4 for the term $\alpha_3$ means that the contribution to the profit $\Psi$ from the third component is 4 when the component is at level 1, 8 when the component is at level 2 and so on. The "interaction" term $2\alpha_3\alpha_4$ implies that profit is further boosted when components 3 and 4 are simultaneously at high performance levels.

For the profit function (11) defined on the state space $\mathcal{D} = \{0, 1, 2, 3\}^4$, if the cutoff point is set at at $c = 28$ we obtain nine minimal nonfailure configurations,

$$\bar{\mathcal{F}}_{28}^* = \{(3, 2, 3, 1), (2, 3, 3, 1), (2, 0, 2, 2), (1, 1, 2, 2),$$
$$(0, 2, 2, 2), (3, 0, 1, 3), (2, 1, 1, 3), (1, 2, 1, 3), (0, 3, 1, 3)\}.$$

Since the profit function is increasing in all the components, any configuration $\beta$ with $\beta \succeq \alpha$, with $\alpha$ in $\bar{\mathcal{F}}_{28}^*$, still belongs to the nonfailure set $\bar{\mathcal{F}}_{28} = \{\alpha \in \mathcal{D} \mid \Psi(\alpha) \geq 28\}$.

To obtain the reliability $\mathcal{R} = \text{Prob}(\bar{\mathcal{F}}_{28})$ and the corresponding bounds we use the improved inclusion–exclusion formula via the Scarf complex, as described above. Since the genericity condition does not hold for $\bar{\mathcal{F}}_{28}^*$ we need to proceed via deformation and ranking:

$$\tilde{\mathcal{F}}_{28}^* = \{(7, 4, 7, 0), (4, 7, 8, 1), (5, 0, 4, 2), (2, 2, 5, 3),$$
$$(0, 5, 6, 4), (8, 1, 0, 5), (6, 3, 1, 6), (3, 6, 2, 7), (1, 8, 3, 8)\}.$$

After labeling the minimal nonfailure points in $\tilde{\mathcal{F}}_{28}^*$ from 1 to 9 the facets of the corresponding Scarf complex can be obtained:

$$\Delta_{28} = \{\{123\}, \{489\}, \{459\}, \{234\}, \{378\}, \{348\}, \{367\},$$
$$\{49\}, \{59\}, \{23\}, \{24\}, \{38\}, \{89\}, \{78\}, \{48\}, \{36\},$$
$$\{67\}, \{45\}, \{34\}, \{13\}, \{37\}, \{12\}, \{9\}, \{8\}, \{7\}, \{6\}, \{5\}, \{4\}, \{3\}, \{2\}, \{1\}\}.$$

The geometrical representation of the Scarf complex is given in Figure 2. Notice that the inclusion–exclusion formula now includes only 31 terms compared to the $2^9 - 1 = 511$ terms that would appear in the complete classical formula.

As an example of the calculation of the reliability formula (8) or the tube bounds (9), we obtain the probability of the orthant associated with the index set $\{489\}$ in the Scarf complex $\Delta_{28}$.

Using the notation introduced in [9], we indicate with $X_i$ the random state of component $i$ and with $\bar{P}_i(j)$ the probability that component $i$ is in a state level "better" than $j$: $\bar{P}_i(j) = \text{Prob}(X_i \geq j)$. In addition we indicate with $p_{ij}$ the probability that component $i$ is exactly in state $j$: $p_{ij} = \text{Prob}(X_i = j)$, $i = 1, \ldots, 4$, $j = 0, \ldots, 3$. Therefore we obtain the following expression for $\bar{P}_i(j)$:

$$\bar{P}_i(j) = \sum_{k=j}^{n_j} p_{ik},$$

where $n_j$ is the maximum performance level of the $j$th component. The index set $I = \{489\}$ in the Scarf complex corresponds to the least common multiple of the three monomials in the (original, before perturbation) minimal failure sets $\bar{F}_{28}^*$ labeled $4, 8$ and $9$, respectively:

$$m_I = \text{lcm}(x_1^1 x_2^1 x_3^2 x_4^2, \ x_1^1 x_2^2 x_3^1 x_4^3, \ x_2^3 x_3^1 x_4^3) = x_1^1 x_2^3 x_3^2 x_4^3.$$

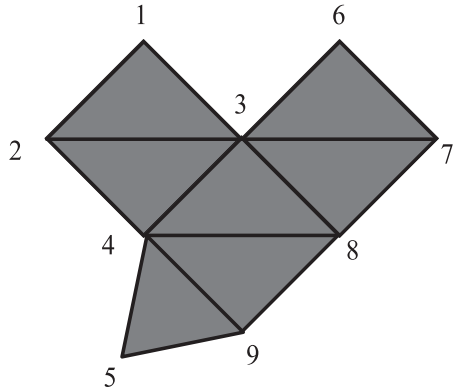

FIG. 2. *Geometric realization of the Scarf complex in the multistate example.*

Thus, in the case of mutual independence among the components, the probability of the orthant associated with the index set $I = \{489\}$ [and the corresponding monomial $m_I = x_1^1 x_2^3 x_3^2 x_4^3$, $\alpha_I = (1, 3, 2, 3)$] is obtained by

$$\mathrm{Prob}(Q_{\alpha_I}) = \prod_{i=1}^{4} \mathrm{Prob}(X_i \geq \alpha_{Ii})$$

$$= \bar{P}_1(1)\bar{P}_2(3)\bar{P}_3(2)\bar{P}_4(3) = p_{23} \cdot p_{43} \cdot \sum_{j=1}^{3} p_{1j} \cdot \sum_{j=2}^{3} p_{3j}.$$

5.5. *The continuous distribution case.* A natural generalization of the multistate system is to components with continuous states. With care, the continuous case can be mapped into the discrete case in such a way that the inclusion–exclusion lemma and associated tube bounds can be inferred back from the discrete to the continuous case.

A simple way to define the mapping is to replace continuous values by their rank along each dimension. In addition one needs to "quantize" the continuous distributions.

Thus let $Z = (Z_1, \ldots, Z_n)$ be an $n$-vector random variable with distribution function

$$F_Z(z) = \mathrm{Prob}(Z \leq z).$$

Define critical points

(12) $$z^{(1)}, \ldots, z^{(m)}$$

and the corresponding orthants

$$Q_{z^{(i)}} = \{z | z > z^{(i)}\}.$$

Thus the nonfailure event is given by

$$\bar{\mathcal{F}} = \bigcup_{i=1}^{m} Q_{z^{(i)}}.$$

Now consider for a particular dimension $i$ the set of all values on that dimension,

$$z_i^{(1)}, \ldots, z_i^{(m)}.$$

Assume the continuous genericity (general position) condition,

(13) $$\text{for all} \quad j \neq j' \implies z_i^{(j)} \neq z_i^{(j')} \quad j = 1, \ldots, m.$$

The condition is not too strong in this case because for continuous variables the critical points (12) can be easily chosen so that condition (13) is satisfied.

SCARF COMPLEX FOR COHERENT SYSTEMS 23

Fix the dimension $i$ and relabel the values so that

$$z_i^{(1)} < z_i^{(2)} < \cdots < z_i^{(m)}.$$

Finally, create $m+1$ states for the discrete variable $X_i$ so that

$$X_i \text{ is in state } j \quad \Longleftrightarrow \quad z_i^{(j)} < Z_i \leq z_i^{(j+1)}.$$

This induces a distribution on the $(m+1)^n$ grid with probabilities

$$\text{Prob} \bigcap_{i=1}^{n} \{z_i^{(j)} < Z_i \leq z_i^{(j+1)}\},$$

which allows the orthant probabilities to be expressed in a straightforward way in terms of the distribution function $F$,

$$\text{Prob}(Q_{z^{(i)}}) = 1 - F(z_1^{(i)}, \ldots, z_n^{(i)}).$$

The Scarf complex and associated tube bounds as obtained in the previous sections can then be derived.

RISK REVIEW DEPARTMENT
FINANCIAL SERVICES AUTHORITY
25 THE NORTH COLONNADE
LONDON E14 5HS
UNITED KINGDOM
E-MAIL: beatrice.giglio@virgin.net

DEPARTMENT OF STATISTICS
LONDON SCHOOL OF ECONOMICS
HOUGHTON STREET
LONDON WC2A 2AE
UNITED KINGDOM
E-MAIL: h.wynn@lse.ac.uk